\documentclass{amsart}
\usepackage[english]{babel}
\usepackage{amsfonts, amsmath, amsthm, amssymb,amscd,indentfirst}

\usepackage{esint}
\usepackage{graphicx}
\usepackage{epstopdf}
\usepackage{amsmath,amssymb,latexsym,indentfirst}
\usepackage{times}
\usepackage{palatino}
\newtheorem{theorem}{Theorem}[section]

\newtheorem{lemma}{Lemma}[section]

\newtheorem{corollary}{Corollary}[section]
\newtheorem{observation}{Observation}[section]
\newtheorem{conjecture}{Conjecture}[section]
\newtheorem{question}{Question}[section]

\newcommand{\dsum}{\displaystyle\sum}

\newcommand{\inte}{{\rm int}}

\begin{document}

\title[A characterization of the delaunay surfaces]{A characterization of the delaunay surfaces}

\author{Ezequiel Barbosa}
\address{Universidade Federal de Minas Gerais (UFMG), Departamento de Matem\'{a}tica, Caixa Postal 702, 30123-970, Belo Horizonte, MG, Brazil.}
\email{ezequiel@mat.ufmg.br}
\author{L.C. Silva}
\address{Universidade Federal de Vi\c{c}osa - campus UFV Florestal, Instituto de Ci\^{e}ncias \linebreak Exatas e Tecnol\'{o}gicas, 35690-000, Florestal, MG, Brazil.}
\email{lcs.mat@gmail.com}
%\thanks{The authors were partially supported by CNPq and CAPES/Brazil grants.}

\begin{abstract}
In this paper we use the Alexandrov Reflection Method to obtain a characterization to embedded CMC capillary annulus $\Sigma^2 \subset \mathbb{B}^3$. In especial, but using a new strategy, we present a new characterization to the critical catenoid. Precisely, we show that $\Sigma \subset \mathbb{B}^3$ being an embedded minimal free boundary annulus in $\mathbb{B}^3$ such that $\partial \Sigma$ is invariant under reflection through a coordinates planes, then $\Sigma$ is the critical catenoid.
\end{abstract}

\maketitle

%%%%%%%%%%%%%%%%%%%%%%%%%%%%%%%%%%%%%%%%%%%%%%%%%%%%%%%%%%%%%%%%%%%%%%%%%%%%%%%%%%%%%%%%%%%%%%%%%%%%%%%%

\section{Introduction}\label{intro}

%%%%%%%%%%%%%%%%%%%%%%%%%%%%%%%%%%%%%%%%%%%%%%%%%%%%%%%%%%%%%%%%%%%%%%%%%%%%%%%%%%%%%%%%%%%%%%%%%%%%%%%%

An important and well-know result in geometry, due to the Nitsche \cite{Nitsche}, state that the only minimal disks, free boundary in $\mathbb{B}^3$, are the equatorial flat discs. In the sense, in \cite{Fraser_Li} we have the 

\begin{conjecture}[Fraser and Li]  \label{Fraser Conjecture}
The critical catenoid is the unique properly embedded free boundary minimal annulus in $\mathbb{B}^3$, up to rotations.
\end{conjecture}

The question above is analogous to the question answered by Nitsche, the difference between the two questions is found in the topology of the surfaces - annular or disc. Philosophically, there exist a parallel between the conjecture (\ref{Fraser Conjecture}) and 

\begin{conjecture}[Lawson]  \label{Lawson Conjecutre}
The Clifford Torus is the only embedded minimal torus in $\mathbb{S}^3$, up to rotations.
\end{conjecture}

The Lawson's conjecture was resolved definitively by Brendle in \cite{Brendle}. However, there was previously a partial proof due to Ros, in \cite{Ros}:

\begin{theorem}[Ros]
Let $\Sigma \subset \mathbb{S}^3$ be an embedded minimal torus, symmetric with respect to
the coordinate hyperplanes of $\mathbb{R}^4$. Then $\Sigma$ is the Clifford torus.
\end{theorem}

In the case of the conjecture (\ref{Fraser Conjecture}), there is a result analogous to that obtained by Ros, due to McGrath \cite{McGrath}:

\begin{theorem}[McGrath]  \label{Teo_McGrath}
Let $\Sigma \subset \mathbb{B}^n$, $n \geq 3$, be an embedded free boundary \linebreak minimal annulus. If $\Sigma$ is invariant under reflection through three hyperplanes, orthogonal to each other, $\Pi_i$, $i = 1, 2, 3$, then $\Sigma$ is the critical catenoid, up to rotation.
%$\partial\Sigma \cap \left( \mathbb{B}^n \setminus \cup_{i = 1}^{3} \Pi_i \right) \neq \emptyset$, 
\end{theorem}

\indent Is know well that, if $\Sigma$ is a minimal surface free boundary in $\mathbb{B}^3$, then it \linebreak coordinate functions are solutions for the Steklov Problem

\begin{align}
\left\{ \begin{array}{cl}
\Delta u = 0, & \textrm{ on } \ \Sigma. \\
\dfrac{\partial u}{\partial \eta} = u, & \textrm{ along } \ \partial\Sigma.
\end{array} \right. 
\end{align}

\indent In him proof, McGrath use the result below, that can be found in \cite{Fraser_Schoen}.

\begin{theorem}[Fraser and Schoen] 
Suppose $\Sigma$ is a free boundary minimal annulus in $\mathbb{B}^n$ such that the coordinate functions are first Steklov eigenfunctions. Then $n = 3$ and $\Sigma$ is \linebreak congruent to the critical catenoid.
\end{theorem}

\indent In this paper, we presented, in the case $n=2$, an improvement for the McGrath Theorem, as consequence of the following result:

\vspace{0.25cm}

\noindent \textbf{Theorem \ref{Resultado Principal}} 
Let $\Sigma^2 \subset \mathbb{B}^3$ be an embedded CMC capillary annulus, such that $\partial\Sigma$ is symmetrical with respect to the coordinated planes, then $\Sigma$ is a delaunay surface.

\vspace{0.25cm}

\indent This theorem makes significant contributions in comparison with the results found in the literature. In their hypotheses, we consider cmc capillary surfaces instead of free boundary minimal surfaces. Thus, with our work, we shed light on the path that leads to the answer to the following question:

\begin{question}
An embedded CMC capillary annulus $\Sigma^2 \subset \mathbb{B}^3$ must be a delaunay surface.
\end{question}

\indent In the especial case, where $\Sigma^2 \subset \mathbb{B}^3$ is an embedded free boundary minimal annulus, in Theorem \ref{Resultado Principal}, we have an improvement for McGrath's Theorem:

\vspace{0.25cm}

\noindent \textbf{Corollary \ref{McGrath melhorado}} 
Let $\Sigma^2 \subset \mathbb{B}^3$ be an embedded free boundary minimal annulus. If $\partial \Sigma$ is symmetrical with respect to the coordinated planes, then $\Sigma$ is the critical catenoid.

\vspace{0.25cm}

Compared to McGrath's results, we assume that $\partial\Sigma$ is invariant under reflection through three orthogonal hyperplanes, in contrast, he assumes such a propriety for $\Sigma$. Thus, when $n=2$, the following corollary is an improvement that we give to McGrath's theorem, in addition to using another strategy, namely, the Alexandrov Reflection Method (ARM), which we'll talk more about in the next section. With this same methodology, we also present a new version, in the embedded case, of the proof from following result, that can be found in \cite{Pyo}, due to Juncheo Pyo.

\vspace{0.25cm}

\noindent \textbf{Theorem \ref{Pyo_diferente} [Pyo]}
Let $\Sigma^2$ be an embedded minimal surface in $\mathbb{R}^3$ with two boundary components and let $\Gamma$ be one component of $\partial\Sigma$. If $\Gamma$ is a circle and $\Sigma$ meets a plane along $\Gamma$ at a constant angle, then $\Sigma$ is part of the catenoid.

\vspace{0.25cm}

%%%%%%%%%%%%%%%%%%%%%%%%%%%%%%%%%%%%%%%%%%%%%%%%%%%%%%%%%%%%%%%%%%%%%%%%%%%%%%%%%%%%%%%%%%%%%%%%%%%%%%%%

\section{Maximum Principles}

%%%%%%%%%%%%%%%%%%%%%%%%%%%%%%%%%%%%%%%%%%%%%%%%%%%%%%%%%%%%%%%%%%%%%%%%%%%%%%%%%%%%%%%%%%%%%%%%%%%%%%%%

\indent Let $A \subset \mathbb{R}^n$ an open set and
\begin{align}  \label{Operador eliptico}
L(w) = \displaystyle\sum_{i,j} a_{ij}(x) w_{ij} + \displaystyle\sum_{i} b_{i}(x) w_{i} + c(x)w
\end{align}

\noindent where $w_{i} := \dfrac{\partial w}{\partial x_i}$, $w_{ij} := \dfrac{\partial w}{\partial x_i \partial x_j}$ and the functions $a_{ij}, b_{i}$ and $c$ are continuous on $\bar{A}$, a \textit{differential elliptic operator on} $A$, i.e., the matrix $[a_{ij}(x)]$ is positive definite for all $x \in A$, that is, 
\begin{align}
0 < \dsum_{i,j = 1}^{n} a_{ij} \xi_i \xi_j , \ \forall \ x \in A, \ \forall \ \xi \in \mathbb{R}^n \setminus \{ 0 \}.
\end{align}

\noindent We called $L$ \textit{uniformly elliptic on} $A$ if, there exist a constant $\kappa$ such that 
\begin{align}
\kappa |\xi|^2 \leq \dsum_{i,j = 1}^{n} a_{ij} \xi_i \xi_j , \ \forall \ x \in A, \ \forall \ \xi \in \mathbb{R}^n \setminus \{ 0 \}.
\end{align}

\indent Now, we will present three maximum principles that we use during this work, especially in one step of the ARM, and can be found in \cite{Wente}. The first of them, for points $x \in \inte A$:

\begin{lemma} \label{PM}
Let $L$ be an elliptic operator as in (\ref{Operador eliptico}) and $w \in C^2(A)$ a function such that
\begin{align}  
L(w) \geq 0, \ \textrm{ on } A. 
\end{align}

\noindent If exist $x_0 \in A$ such that $w(x_0) = 0$ and $w \leq 0$ on $A$, then $w \equiv 0$ on $A$.
\end{lemma}

\indent The second, for points $x \in \partial A$ such that $\partial A$ is of class $C^1$.

\begin{lemma} \label{PM2}
Let $L$ be an uniformly elliptic operator as in (\ref{Operador eliptico}), let $A$ be a region in $\mathbb{R}^2$ and suppose that in a neighborhood of $x_0 \in \partial A$, $\partial A$ is of class $C^1$. If 
\begin{align}
L(w) \geq 0, \textrm{ on } A, 
\end{align}

\noindent $w(x_0) = 0$, $w(x) \leq 0, \ \forall \ x \in \bar{A}$, and $\frac{\partial w}{\partial \nu} = 0$, where $\nu$ is the inward normal derivative, then $w \equiv 0$, on $\bar{A}$.
\end{lemma}

\indent Finally, the third, for points $x \in \partial A$, at a corner.

\begin{lemma}[Serrin's Boundary Point Lemma at a Corner \cite{Serrin}] \label{PM3}
Let $A \subset \mathbb{R}^2$ be a bounded region which has a $C^2$ boundary in a neighborhood of $x_0 \in \partial A$. Consider $T$ be a normal plane to $\partial A$ at $x_0$ and $A^{+}$ be that component of $A$ lying on one side of $T$ which contains $x_0$ in its closure. Let $L$ be an uniformly elliptic operator on $A^{+}$. Suppose also that 
\begin{align}
| \displaystyle\sum_{i,j} a_{ij}(x) \xi_{i} \nu_{j} | \leq K \cdot [ |(\xi \cdot \nu )| + |\xi|d ] 
\end{align}

\noindent for some constant $K > 0$, all $x \in \overline{A}^{+}$, any $\xi = (\xi_1, ..., \xi_n)$, where $\nu = (\nu_1, ..., \nu_n)$ is an unit normal to $T$, and where $d$ is the distance from $x$ to $T$.

\indent Let $w \in C^{2}(\overline{A}^{+})$ satisfy $L(w) \geq 0$ on $\overline{A}^{+}$ and suppose that $w(x_0) = 0$, $w(x) \leq 0$, for all $x \in \overline{A}^{+}$, and that $\partial w / \partial s = \partial^{2} w / \partial s^{2} = 0$, in any direction which enters $A^{+}$ non-tangentially at $x_0$. 
\end{lemma}

\indent Now, as an application of maximum principles above, we presented the steps of ARM.

\begin{enumerate}

\item Consider a subsidiary plane $P$ and an arbitrary family, $P_{\lambda}$, $\lambda \in \mathbb{R}$, of parallel planes with each other, and orthogonal to $P$.

\item Varying the parameter $\lambda$, a \textit{moving planes} process is started by means of family $P_{\lambda}$. For some $\lambda \in \mathbb{R}$, $P_{\lambda} \cap \Sigma \neq \emptyset$ and can be considered the reflection, through $P_{\lambda}$, of the part of $\Sigma$ surpassed by $P_{\lambda}$.

\item For a critical parameter, $\lambda^{*}$, it is considered the reflection through $P_{\lambda^{*}}$ of the part of $\Sigma$ surpassed by $P_{\lambda^{*}}$, see Figure \ref{steps_ARM}.

\begin{figure}[ht]  
\centering
\includegraphics[scale=0.2]{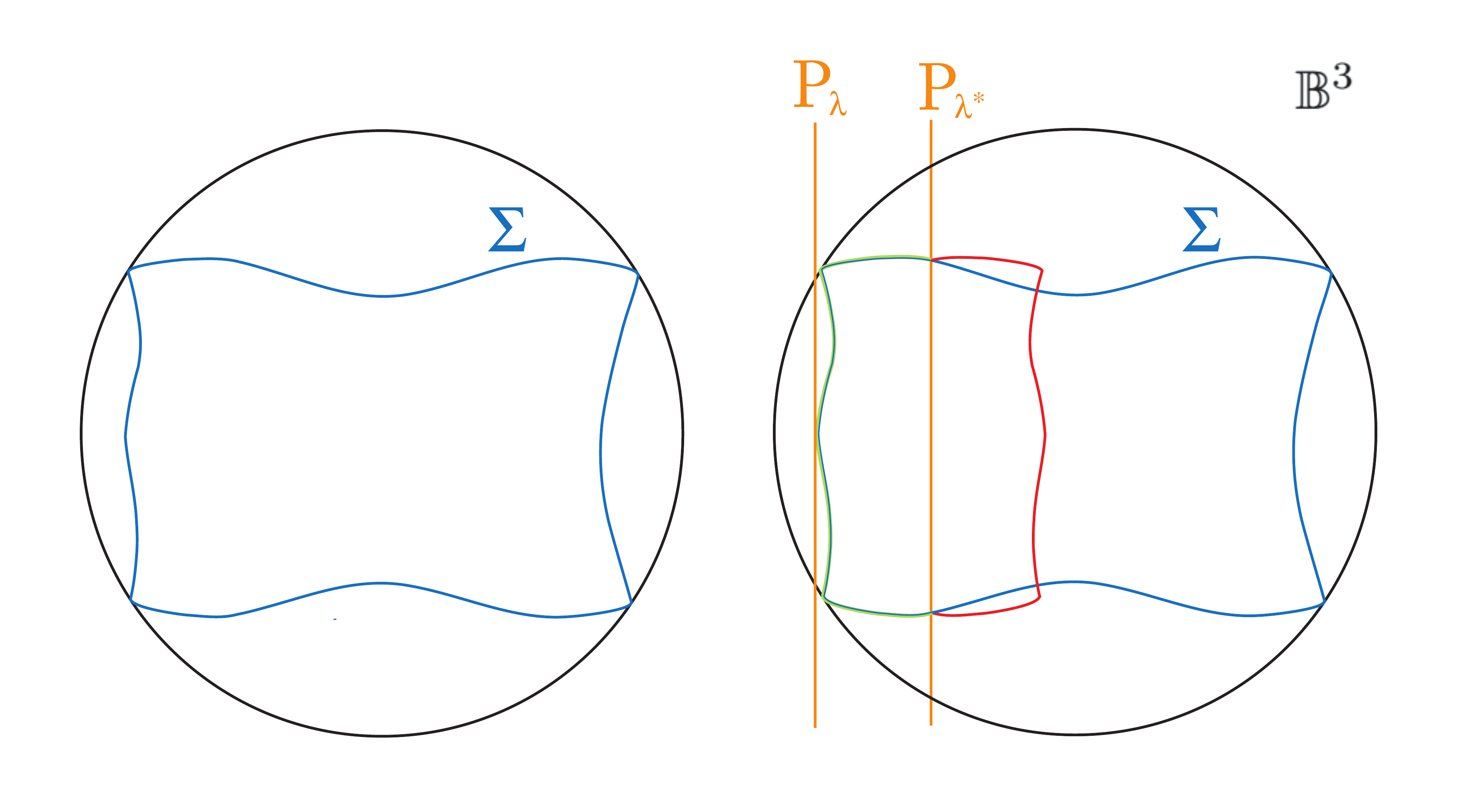}
\caption{In green, the part of $\Sigma$ surpassed by $P_{\lambda^{*}}$. And, in red, the reflection, through $P_{\lambda^{*}}$, of the part of $\Sigma$ surpassed by $P_{\lambda^{*}}$.}
\label{steps_ARM}
\end{figure}

\item Considering an appropriate coordinated system, we use an appropriate maximum principle and it is concluded that the reflection, through $P_{\lambda^{*}}$, of the part of $\Sigma$ surpassed by $P_{\lambda^{*}}$ coincide, locally, with the part of $\Sigma$ non surpassed by $P_{\lambda^{*}}$.

\item The single continuation principle is used and it is concluded that the reflection, through $P_{\lambda^{*}}$, of the part of $\Sigma$ surpassed by $P_{\lambda^{*}}$ coincide with the part of $\Sigma$ non surpassed by $P_{\lambda^{*}}$.

\item Finally, from arbitrariness of $P_{\lambda}$, it is concluded that $\Sigma$ is symmetrical rotationally.

\end{enumerate}

\indent For more examples of this method see \cite{Wente} and \cite{Koiso}. A natural question around the steps above:

\begin{question}  \label{question_moving planes}
How to determine the critical parameter $\lambda^{*}$?
\end{question}

\indent Consider 
\begin{align} \label{Lambda}
\Lambda \textrm{ the region bounded by }  C_{+} \cup \Sigma \cup C_{-} \ \subset \ \mathbb{B}^{3},
\end{align}

\begin{figure}[ht]  
\centering
\includegraphics[scale=0.3]{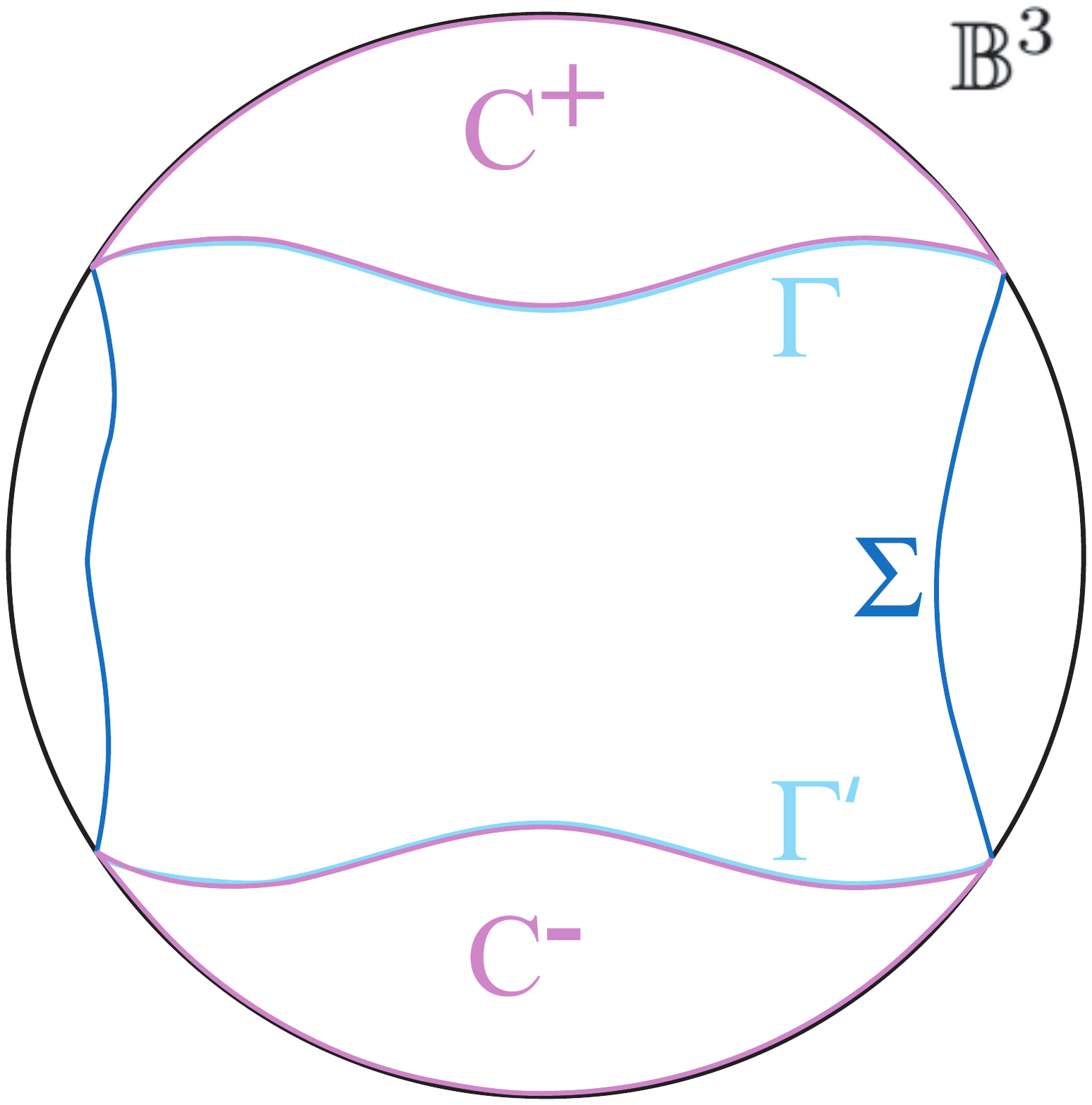}
\caption{$\Lambda$ is the connected region bounded by $C_{+} \cup \Sigma \cup C_{-} \subset \mathbb{B}^{3}$.}
\label{calotas}
\end{figure}
\noindent where $C^{+}$ is the upper portion of $\mathbb{S}^2$ such that $\partial C^{+} = \Gamma$; $C^{-}$ is the lower portion of $\mathbb{S}^2$ such that $\partial C^{-} = \Gamma'$. As $\Sigma$ is embedded, $\Lambda$ is connected (Figure \ref{calotas}).  

%so let $\Pi \cap \Gamma = \{ p, q \}$. Therefore, $< p , e_3 > = < q , e_3 >$, since $\partial\Sigma$ is $G$-invariant, where $\{ e_1, e_2, e_3 \}$ is the canonic base in $\mathbb{R}^3$. 
%  $\Pi_i$, $i \in \{ 1, 2 \}$, e orthogonal to $\Pi_3$. We have to $\Sigma'_{+}(T_\lambda)$ will be contained in $\Lambda$ until one of the four situations below happens.
%Consider $P_{\lambda}$ a family of parallel planes to $P_0$.

\indent Define
\begin{align}
\lambda^{-} = \min \, \{ \lambda \in \mathbb{R} \ ; \ P_{\lambda} \cap \Sigma \neq \emptyset \},
\end{align}
\begin{align}
\lambda^{+} = \max \ \{ \lambda \in \mathbb{R} \ ; \ P_{\lambda} \cap \Sigma \neq \emptyset \},
\end{align}

\noindent and observe that, as $\Sigma \subset \mathbb{B}^3$, so $-1 \leq \lambda^{-} < \lambda^{+} \leq 1$. To better organize the text, consider the following definition: 

\begin{description}

\item[(i)] $\Sigma_{\lambda}$ being \textit{the part of} $\Sigma$ \textit{between} $P_{\lambda^{-}}$ \textit{and} $P_{\lambda}$, $\lambda \in (\lambda^{-}, \lambda^{+})$, is that, the part of $\Sigma$ surpassed by $P_{\lambda}$; 

\item[(ii)] $\widetilde{\Sigma}_{\lambda}$ being \textit{the reflection of} $\Sigma_{\lambda}$ \textit{through} $P_{\lambda}$;

\item[(iii)] $\Sigma\!\setminus\!\Sigma_{\lambda}$ being \textit{the part of} $\Sigma$ \textit{between} $P_{\lambda}$ \textit{and} $P_{\lambda^{+}}$, $\lambda \in (\lambda^{-}, \lambda^{+})$, is that, the part of $\Sigma$ non surpassed by $P_{\lambda}$.  

\end{description}

%Let $\delta > 0$, small enough. For $\lambda = \lambda^{-} + \delta$, the reflection, through $P_{\lambda}$, of the part of $\Sigma$ surpassed by $P_{\lambda}$ is contained in $\Lambda$. This is true until one of the \textit{touching possibilities} occurs:

\begin{figure}[ht]  
\centering
\includegraphics[scale=0.18]{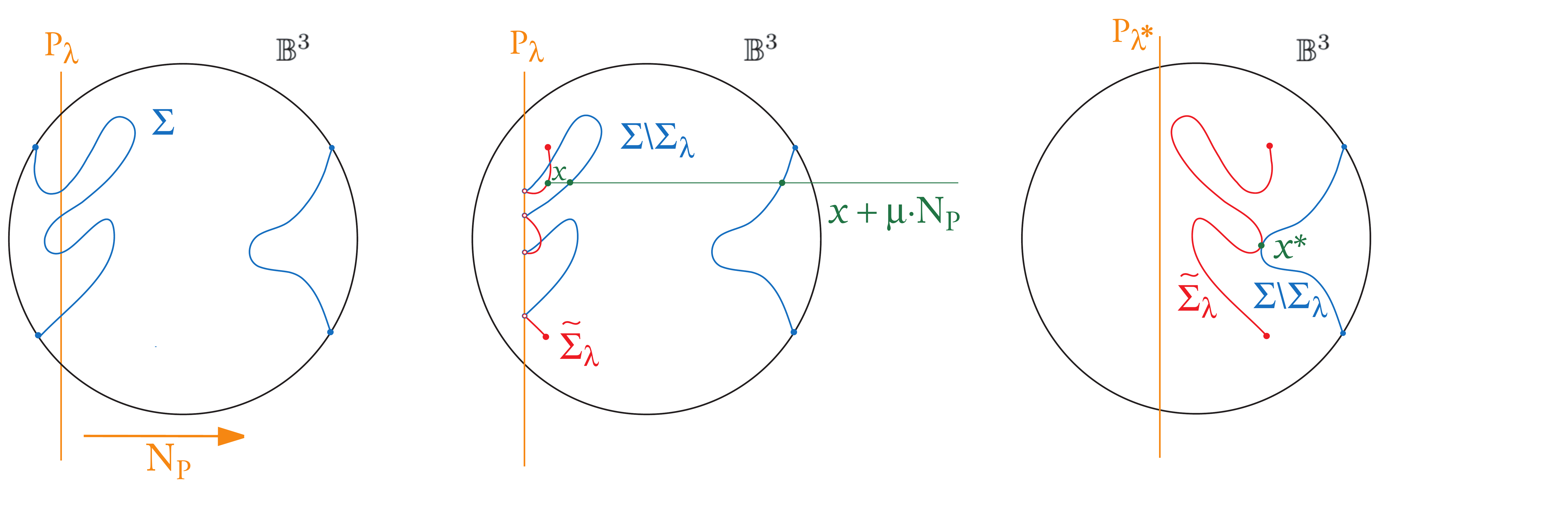}
\caption{The moving planes process in two moments, $\lambda$ and $\lambda^{*}$.}
\label{definitivamente_2}
\end{figure}

For some value of parameter $\lambda$, called $\lambda^{*}$, we say that the reflected part, $\widetilde{\Sigma}_{\lambda}$, \textit{definitely extrapolates} $\Lambda$ if, 
\begin{align}
\exists \ x^{*} \in  \widetilde{\Sigma}_{\lambda} \ ; \ x^{*} + \mu \cdot N_{P} \not\in \Lambda, \ \forall \ \mu > 0,
\end{align}

\noindent where $N_{P}$ is the unit normal vector to family $P_{\lambda}$, pointing in the sense of increasing $\lambda$. Thus, it is defined the critical parameter of moving planes process with respect to family $P_{\lambda}$, see Figure \ref{definitivamente_2}. There are the following possibilities for this extrapolation:

\begin{description}

\item[(P1)] At a point $x^{*} \ \textrm{on} \ \inte(\widetilde{\Sigma}_{\lambda^{*}}) \!\cap \inte(\Sigma\!\setminus\!\Sigma_{\lambda^{*}})$.

\item[(P2)] At a point $x^{*} \in \partial \widetilde{\Sigma}_{\lambda^{*}} \!\cap \partial (\Sigma\!\setminus\!\Sigma_{\lambda^{*}})$.

\item[(P3)] At a point $x^{*}$ such that $T_{x^{*}}\Sigma \perp P_{\lambda^{*}}$.

\item[(P4)] At a point $x^{*}$ such that $T_{p}\partial\Sigma \perp P_{\lambda^{*}}$.

\end{description}

Note that, we should not worry with the possibility of $\widetilde{\Sigma}_{\theta, \lambda}$ definitively extrapolates $\Lambda$ by $C^{+}$ or $C^{-}$ and also in the possibility of a point along $\partial\widetilde{\Sigma}_{\theta, \lambda}$ definitively extrapolates $\Lambda$, at a point $p \in \inte \Sigma\!\setminus\!\Sigma_{\lambda}$, because $\partial\Sigma$ is invariant under reflection through of the three coordinated hyperplanes and due to spherical geometry (this is another relevance of the $\partial\Sigma$ symmetry). 

\begin{observation} We created and adopted the concept \textit{definitively extrapolates} instead of \textit{touching}, the latter already existing in the literature, to avoid the possibility of the touch occurring at the intersection of a 
boundary point of $\widetilde{\Sigma}_{\lambda}$ with an interior point of $\Sigma\!\setminus\!\Sigma_{\lambda}$. 
\end{observation}

\indent Once the Question \ref{question_moving planes} has been answered, we can present our results.

%%%%%%%%%%%%%%%%%%%%%%%%%%%%%%%%%%%%%%%%%%%%%%%%%%%%%%%%%%%%%%%%%%%%%%%%%%%%%%%%%%%%%%%%%%%%%%%%%%%%%%%%

\section{The Proof of the Theorem \ref{Resultado Principal}}

%%%%%%%%%%%%%%%%%%%%%%%%%%%%%%%%%%%%%%%%%%%%%%%%%%%%%%%%%%%%%%%%%%%%%%%%%%%%%%%%%%%%%%%%%%%%%%%%%%%%%%%%

\indent Let $\Sigma \subset \mathbb{B}^3$ an embedded free boundary minimal annulus such that 

\begin{align}  \label{admissibilidade}
\inte(\Sigma) \subset \inte(\mathbb{B}^3) \\
\partial\Sigma = \Gamma \cup \Gamma' \subset \partial\mathbb{B}^3
\end{align}

\noindent where $\Gamma$ and $\Gamma'$ are the connected components of the boundary of $\Sigma$. Consider a hyperplane $\Pi \subset \mathbb{R}^3$ and $R_{\Pi}$ the map such that $R_{\Pi}(x)$ is the orthogonal reflection of $x$ through $\Pi$. If $R_{\Pi} (\Sigma) = \Sigma$, we say that $\Sigma$ is $\Pi \textit{-invariant}$. Note that, the map $R_{\Pi} : \Sigma \rightarrow \Sigma$ is an isometry such that $\partial\Sigma \mapsto \partial\Sigma$ and $\inte(\Sigma) \mapsto \inte(\Sigma)$. From now on, consider 

\begin{align}
\Pi_i := \left\{ (x_1, x_2, x_3) \ | \ x_i = 0 \right\},
\end{align}

%\noindent and $j_{*} = 3$, i.e., $\Gamma' = R_{\Pi_{3}}(\Gamma)$.

\indent Let $G = \{ R_{\Pi_1}, R_{\Pi_2}, R_{\Pi_3} \}$ be the group of the reflection with respect to the \linebreak coordinate planes. We say that $\Sigma$ is $G\textit{-invariant}$ if, $R_{\Pi_i} (\Sigma) = \Sigma$, for all $i \in \{ 1,2,3 \}$. In this paper, we consider $\partial \Sigma = \Gamma \cup \Gamma'$ $G$-invariant.

In this moment, we prove that the $G$-invariant propriety of $\partial\Sigma$ implies that it intersects the interior of each of the eight octants and there exist a plane such that $\Gamma'$ is the reflection of $\Gamma$ through it.

\begin{lemma}  \label{meu lema}
Let $\Sigma^2 \subset \mathbb{B}^3$ an embedded annulus such that $\partial \Sigma$ is $G$-invariant. Then,

\begin{description}

\item[(i)] there exist $i,j \in \{ 1,2,3 \}$, $i \neq j$, such that  
\begin{align}  \label{meu lema 1}
\Gamma = R_{\Pi_{i}}(\Gamma) = R_{\Pi_{j}}(\Gamma), \\
\Gamma' = R_{\Pi_{i}}(\Gamma') = R_{\Pi_{j}}(\Gamma'),
\end{align}

\noindent and 

\item[(ii)] there exists $k \in \{ 1,2,3 \}$, $k \notin \{ i, j \}$, such that 
\begin{align}  \label{meu lema 2}
\Gamma' = R_{\Pi_{k}}(\Gamma) \ \textrm{ and } \ \Gamma = R_{\Pi_{k}}(\Gamma')
\end{align} 

\end{description}
\end{lemma}

\noindent \textbf{Proof of Lemma \ref{meu lema}:} Let $\partial \Sigma = \Gamma \cup \Gamma'$, where $\Gamma$ and $\Gamma'$ are the connected components of the boundary of $\Sigma$, and 
\begin{align}
\Gamma \cap \mathcal{O} =: \gamma : [0,1] \rightarrow \partial\Sigma,
\end{align}

\noindent where $\mathcal{O} = \{ (x_1, x_2, x_3) \in \mathbb{R}^3 ; \ x_1, x_2, x_3 \geq 0 \}$, the part of $\Gamma$ contained in the first octant. As $\partial\Sigma$ is $G$-invariant, we can find it by putting together the possible reflections of $\gamma$, i.e., there are $i,j,k \in \{ 1, 2, 3 \}$, different from each other, such that
\begin{align}  \label{bordo_fundamental}
\partial\Sigma = \tilde{\gamma} \cup R_{\Pi_k}( \tilde{\gamma} ) 
\end{align}

\noindent where 
\begin{align}  \label{gamma_tilde}
\tilde{\gamma} := \gamma \cup R_{\Pi_i}(\gamma) \cup R_{\Pi_j}(\gamma) \cup (R_{\Pi_i} \!\circ\! R_{\Pi_j})(\gamma) \subset \partial\Sigma
\end{align}

\begin{figure}[ht]
\begin{minipage}[c]{0.45 \textwidth}
\includegraphics[scale=0.4]{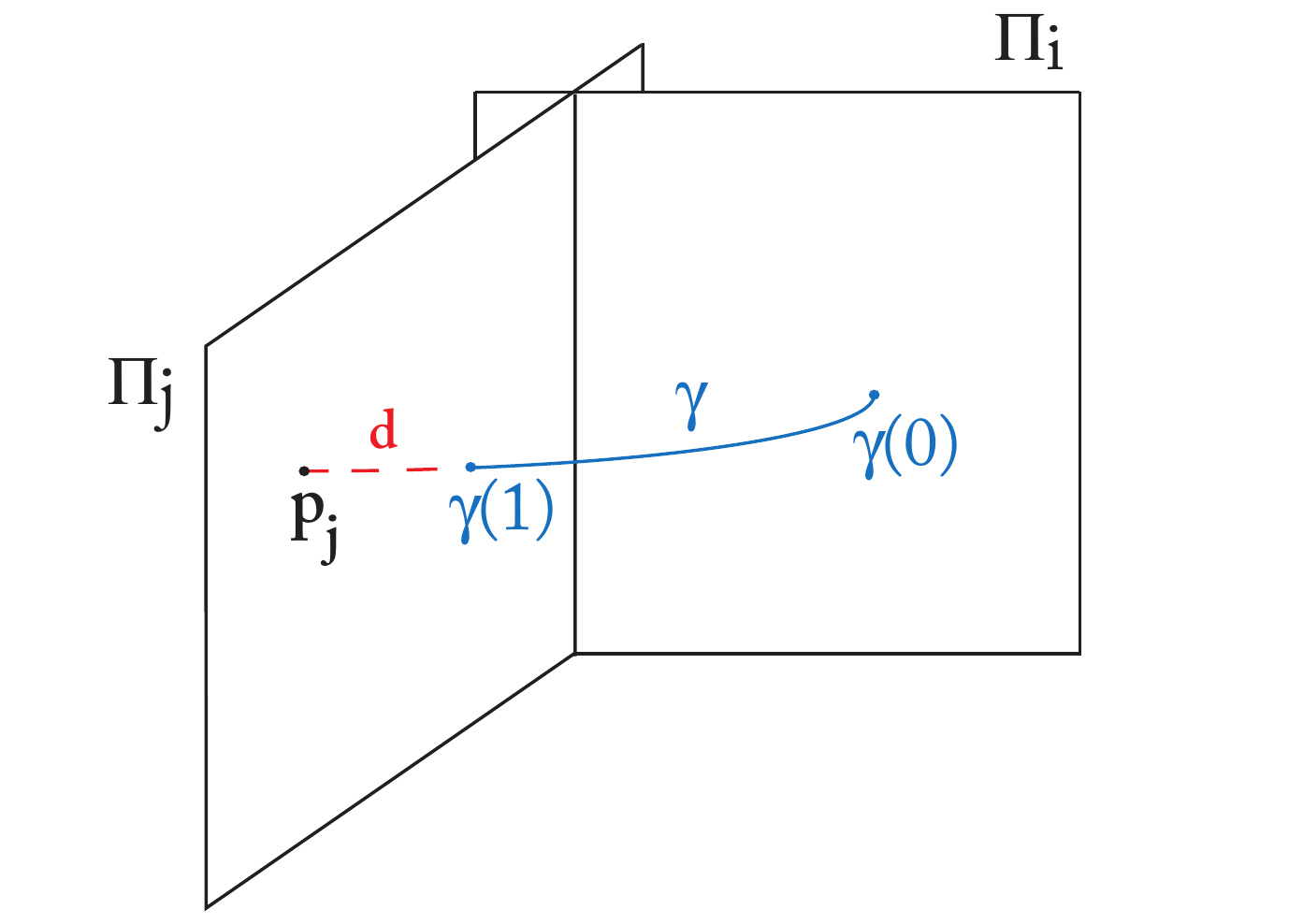}
\caption{$\gamma(1) \not\in \Pi_j$.}
\label{Figura_Lema}
\end{minipage}
\begin{minipage}[c]{0.45 \textwidth}
\includegraphics[scale=0.4]{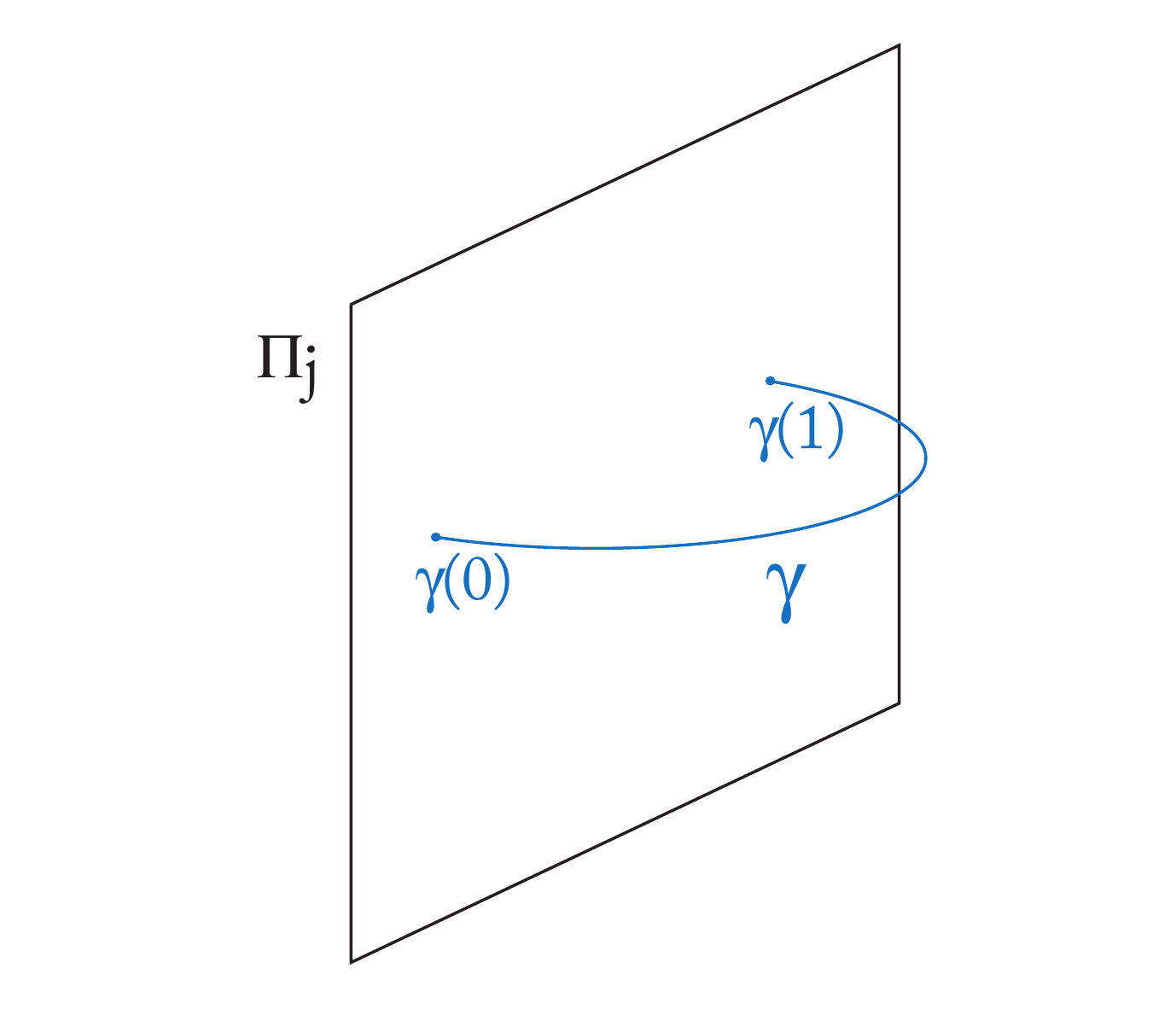}
\caption{$\Pi_i = \Pi_j$ and $\gamma(0), \gamma(1) \in \Pi_j$.}
\label{Figura_Lema_2}
\end{minipage}
\end{figure}

\indent As $\Sigma$ is embedded, $\gamma(0,1)$ doesn't intersect $\Pi_i$, $i \in \{ 1, 2, 3 \}$, since $\partial\Sigma$ is $G$-invariant, otherwise $\Gamma$ would have self intersections. Once $\Gamma$ and $\Gamma'$ are closed curves, $\gamma(0) \in \Pi_i$ and $\gamma(1) \in \Pi_j$, $i \neq j$, because $\partial\Sigma$ is $G$-invariant. Indeed, if $\gamma(1) \not\in \Pi_j$, then there exist $p_j \in \Pi_j$ such that $d(\gamma(1),p_j) = d > 0$ and follows from (\ref{bordo_fundamental}) and (\ref{gamma_tilde}) that, $\partial\Sigma$ would not be the union of closed curves and this would be a contradiction, see Figure \ref{Figura_Lema}. 

In the other hand, if $\gamma(0), \gamma(1) \!\in\! \Pi_i = \Pi_j$, see Figure \ref{Figura_Lema_2}, then the curve $\gamma \cup R_{\Pi_i}(\gamma) := \beta : [a,b] \rightarrow \partial\Sigma$ would be a closed curve contained in $\partial\Sigma$. Thus, the curves $\beta$, $R_{\Pi_j}(\beta)$, $R_{\Pi_k}(\beta)$ and $(R_{\Pi_j} \circ R_{\Pi_k})(\beta)$ would be closed curves contained in $\partial\Sigma$, but this is also a contradiction, since $\Sigma$ is a topological annulus.

\indent Then, if we define $\Gamma := \tilde{\gamma}$, we have $\textbf{(i)}$ and $\textbf{(ii)}$.
\begin{flushright}
$\blacksquare$
\end{flushright}

%\begin{align}
%\Gamma = \gamma \cup R_{\Pi_i}(\gamma) \cup R_{\Pi_j}(\gamma) \cup R_{\Pi_i} \circ R_{\Pi_j}(\gamma) \\
%\Gamma' = R_{\Pi_{k}}(\gamma) ; \ k \in \{ 1, 2, 3 \} \setminus \{ i, j \}.
%\end{align}

\begin{figure}[ht]  
\centering
\includegraphics[scale=0.3]{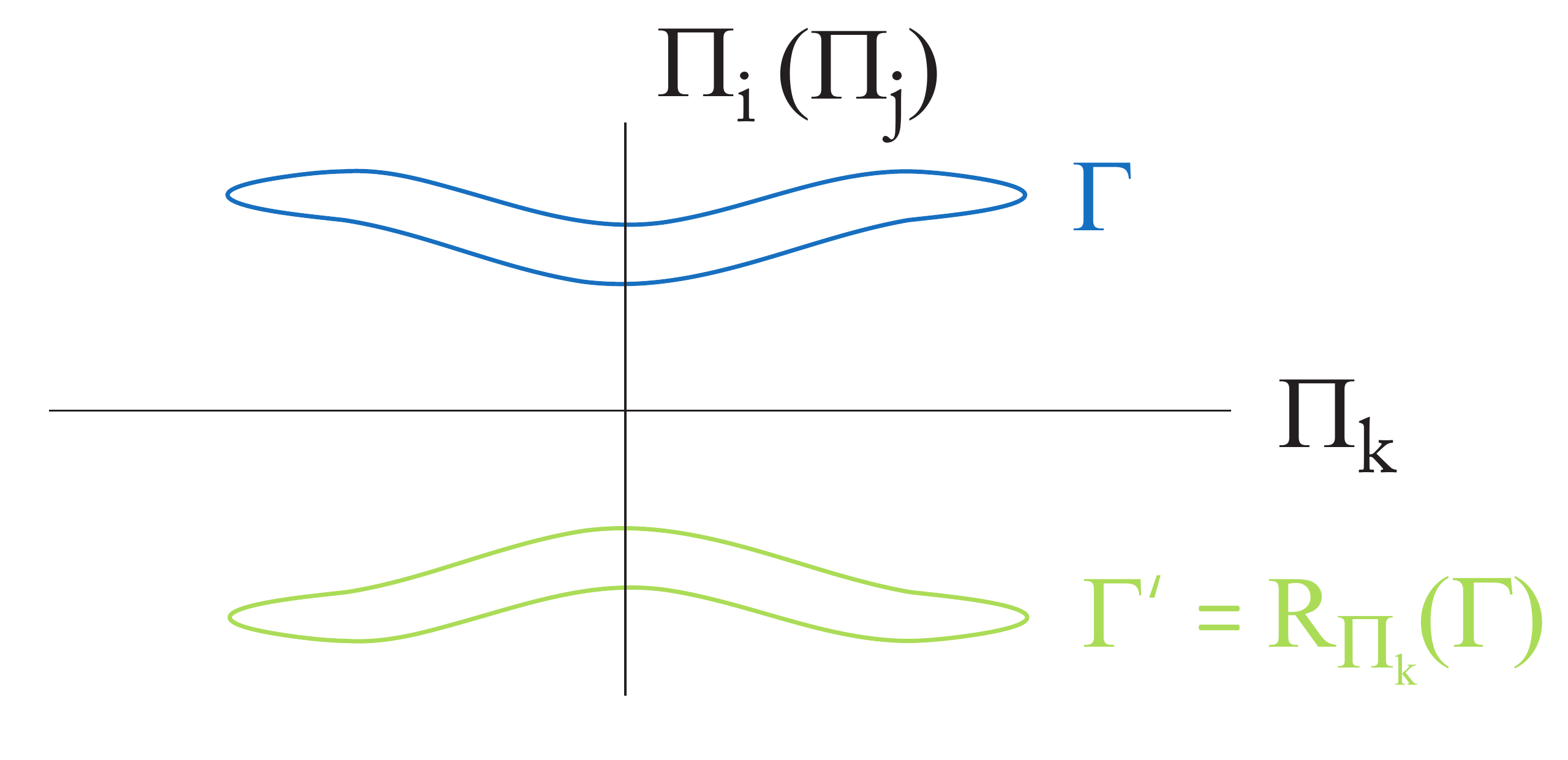}
\caption{$\Gamma' = R_{\Pi_{k}}(\Gamma) ; \ k \in \{ 1, 2, 3 \} \setminus \{ i, j \}$.}
\label{bordo sigma}
\end{figure}
\indent Let $\mathcal{T_{\lambda}}$, $\lambda \in \mathbb{R}$, be a family of planes parallel to each other, where there is a relationship 1-1 between $\lambda \in \mathbb{R}$ and each plane $T \in \mathcal{T}_{\lambda}$. We call \textit{moving planes}, the process of varying the parameter $\lambda$, from the geometric view point, we have a movement between this parallel planes.

\begin{theorem}  \label{Resultado Principal}
Let $\Sigma^2 \subset \mathbb{B}^3$ be an embedded CMC capillary annulus, such that $\partial\Sigma$ is \linebreak symmetrical with respect to the coordinated planes, then $\Sigma$ is a delaunay surface.
\end{theorem}

\noindent \textbf{Proof of Theorem \ref{Resultado Principal}:} Let $\Sigma \subset \mathbb{B}^3$ an embedded CMC capillary annulus such that $\partial \Sigma = \Gamma \cup \Gamma'$ is $G$-invariant and
\begin{align} 
\inte(\Sigma) \subset \inte(\mathbb{B}^3) \\
\partial\Sigma = \Gamma \cup \Gamma' \subset \partial\mathbb{B}^3
\end{align}

\indent As $\partial\Sigma$ is $G$-invariant, follows from Lemma \ref{meu lema} that there exists a coordinated plane, without loss of generality, let's say $\Pi_3$, such that
\begin{align}  \label{bordo espelhado_1} 
\Gamma' = R_{\Pi_3}(\Gamma) \ \ \textrm{and} \ \ \Gamma = R_{\Pi_3}(\Gamma')
\end{align}

\noindent and
\begin{align}  \label{bordo espelhado_2}
\Gamma = R_{\Pi_1}(\Gamma) = R_{\Pi_2}(\Gamma)  \ \ \textrm{and} \ \ \Gamma' = R_{\Pi_1}(\Gamma') = R_{\Pi_2}(\Gamma').
\end{align}

\indent Consider the family $\mathcal{T}_{\theta}$ of orthogonal planes to $\Pi_3$ and parallels to each other, such that
\begin{align}
T_{\theta, \lambda} = \{ (x_1, x_2 , x_3) \in \mathbb{R}^3 ; \cos \theta \cdot x_1 + \sin \theta \cdot x_2 = \lambda \} \in \mathcal{T}_{\theta},
\end{align}  

\noindent where $\theta \in [0, \pi)$ and $\lambda \in \mathbb{R}$. Let
\begin{align}
\lambda_{\theta}^{-} = \min \, \{ \lambda \in \mathbb{R} \ ; \ T_{\theta, \lambda} \cap \Sigma \neq \emptyset \},
\end{align}
\begin{align}
\lambda_{\theta}^{+} = \max \ \{ \lambda \in \mathbb{R} \ ; \ T_{\theta, \lambda} \cap \Sigma \neq \emptyset \},
\end{align}

\indent Define $\Sigma_{\theta, \lambda}$ as the part of $\Sigma$ between $T_{\theta, \lambda_{\theta}^{-}}$ and $T_{\theta, \lambda}$, $\lambda \in (\lambda_{\theta}^{-}, \lambda_{\theta}^{+})$, precisely
\begin{align}
\Sigma_{\theta, \lambda} := \{ x \in \Sigma \ ; \ \lambda_{\theta}^{-} \leq \cos\theta \cdot x_1 + \sin\theta \cdot x_2 \leq \lambda \}.
\end{align}

\begin{figure}[ht]  
\centering
\includegraphics[scale=0.3]{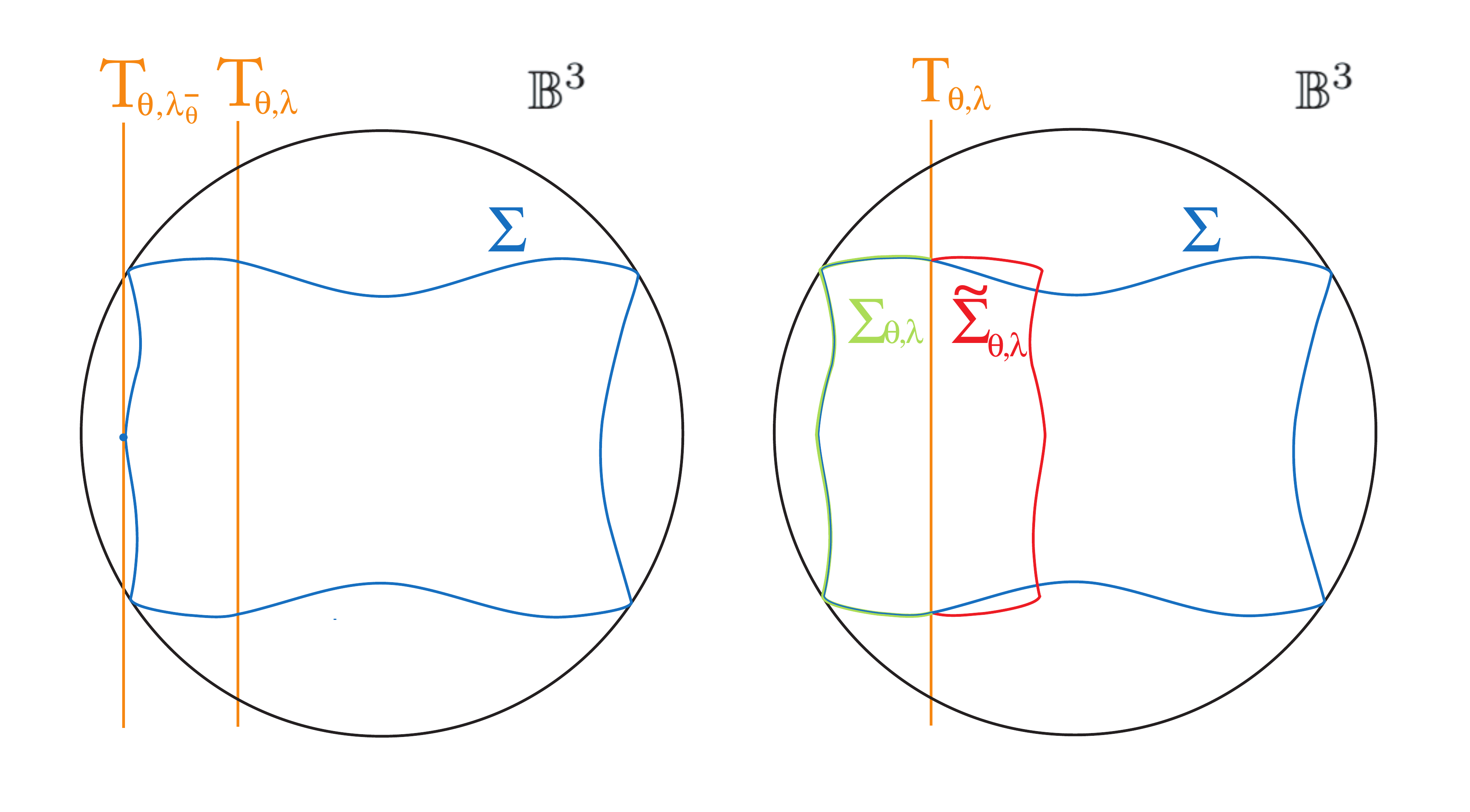}
\caption{$\Sigma_{\theta, \lambda}$, $\widetilde{\Sigma}_{\theta, \lambda}$ and $\Sigma\!\setminus\!\Sigma_{\theta, \lambda}$.}
\label{reflection_Sigma}
\end{figure}

\noindent We will call $\widetilde{\Sigma}_{\theta, \lambda}$ the reflection of $\Sigma_{\theta, \lambda}$ through $T_{\theta, \lambda}$, i.e.,
\begin{align}
\widetilde{\Sigma}_{\theta, \lambda} := \{ x \in \Sigma \ ; \ \lambda \leq \cos\theta \cdot x_1 + \sin\theta \cdot x_2 \leq \lambda + \lambda_{\theta}^{-} \},
\end{align}

\noindent and $\Sigma\!\setminus\!\Sigma_{\theta, \lambda}$ the part of $\Sigma$ between $T_{\theta, \lambda}$ and $T_{\theta, \lambda_{\theta}^{+}}$, $\lambda \in (\lambda_{\theta}^{-}, \lambda_{\theta}^{+})$, see Figure \ref{reflection_Sigma}.  

\begin{figure}[ht]  
\centering
\includegraphics[scale=0.4]{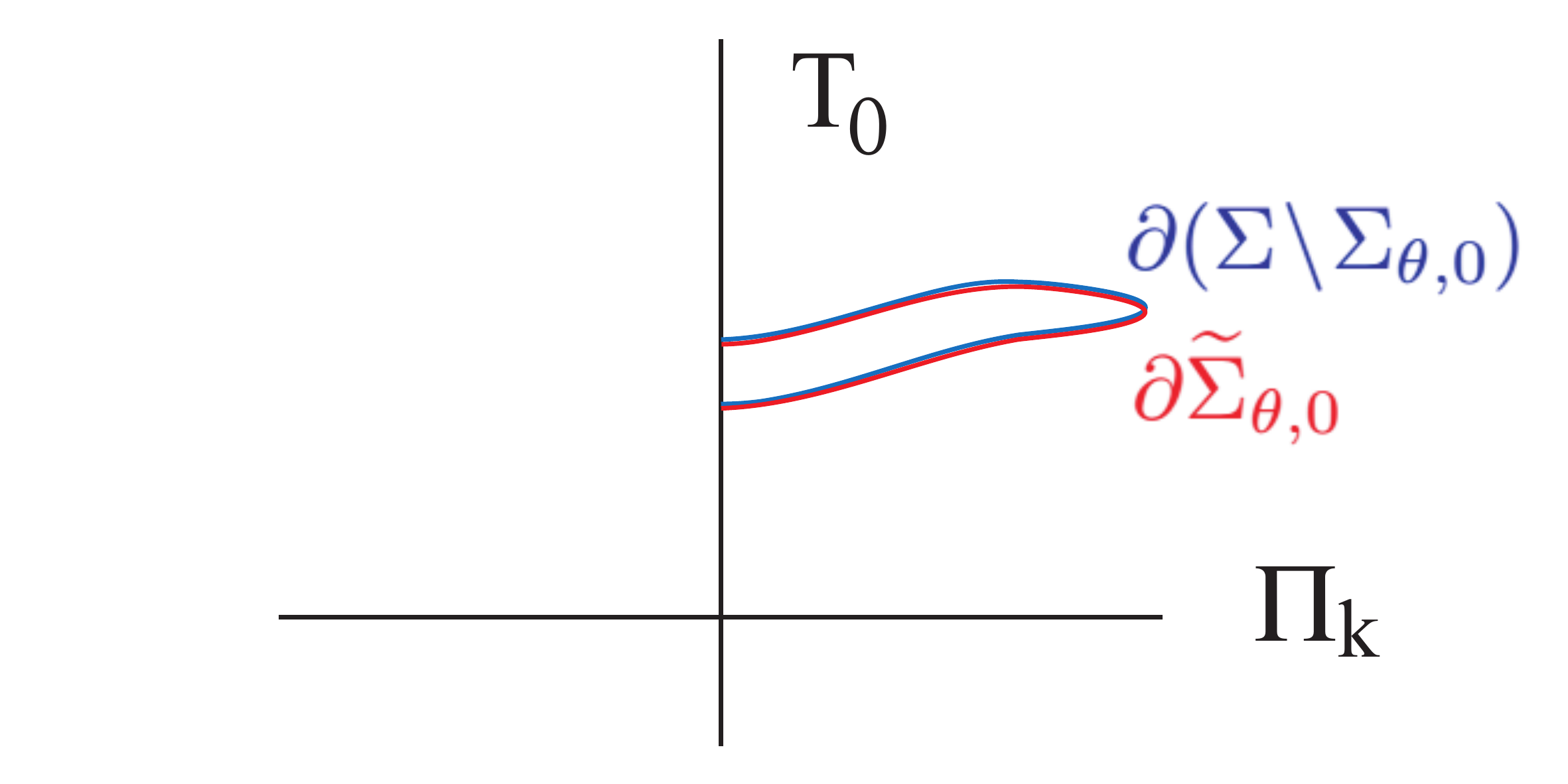}
\caption{$\partial\widetilde{\Sigma}_{\theta, 0} = \partial(\Sigma \! \setminus \! \Sigma_{\theta, 0})$, for $\theta \in \{ 0, \frac{\pi}{2} \}$.}
\label{bordo sigma refletido}
\end{figure}

\indent For example, note that for $\theta = \lambda = 0$, $T_{0,0} = \Pi_1$ and for $\theta = \frac{\pi}{2}$ and $\lambda = 0$, $T_{\frac{\pi}{2},0} = \Pi_2$. In these cases, $\partial \widetilde{\Sigma}_{\theta, \lambda} = \partial \left( \Sigma \setminus \Sigma_{\theta, \lambda} \right)$, see Figure \ref{bordo sigma refletido}. 

\indent Let $x^{*}$ the extrapolation point, for some of the possibilities \textbf{(P1)} $\sim$ \textbf{(P4)}. Consider a coordinate system such that $x^{*}=(0,0,0)$ and smooth functions $u,v : \bar{A} \rightarrow \mathbb{R}$, where $A \subset \mathbb{R}^2$ is an open set and $(0,0) \in \bar{A}$, such that 
\begin{align}
u(0,0) = v(0,0) = 0
\end{align}

\noindent and $\Sigma\!\setminus\!\Sigma_{\theta, \lambda_{\theta}^{*}} = \textrm{graph}(u)$ and $\widetilde{\Sigma}_{\theta, \lambda_{\theta}^{*}} = \textrm{graph}(v)$ in a neighborhood of $(0,0)$. Note that, as $\Sigma\!\setminus\!\Sigma_{\theta, \lambda_{\theta}^{*}}$ and $\widetilde{\Sigma}_{\theta, \lambda_{\theta}^{*}}$ are CMC (for the same constant), $u$ and $v$ satisfy the same CMC equation. Hence, the function $w = v - u$ satisfy a homogeneous linear elliptic pde, see \cite{Wente}.

\indent In the possibility $\textbf{(P1)}$, define a coordinate system such that $T_{x^{*}} \Sigma = \{ z = 0 \}$, where the axis $z$ pointing to $T_{\theta, \lambda_{\theta}^{*}}$ and use the Lemma \ref{PM} to conclude that $w = 0$ in a neighborhood of $(0,0)$, i.e., $\widetilde{\Sigma}_{\theta, \lambda_{\theta}^{*}} = \Sigma\!\setminus\!\Sigma_{\theta, \lambda_{\theta}^{*}}$ in a neighborhood of $x^{*}$. 

\indent In $\textbf{(P2)}$, define a coordinate system such that $T_{x^{*}} \partial\Sigma = \{ x = z = 0 \}$ and \linebreak $T_{x^{*}} \Sigma = \{ z = 0 \} $, where the axis $z$ pointing to $T_{\theta, \lambda_{\theta}^{*}}$ and axis $x$ point to $\inte A$. So, use the Lemma \ref{PM2} to conclude that $w = 0$ in a neighborhood of $(0,0)$, i.e., $\widetilde{\Sigma}_{\theta, \lambda_{\theta}^{*}} = \Sigma\!\setminus\!\Sigma_{\theta, \lambda_{\theta}^{*}}$ in a neighborhood of $x^{*}$.

\indent For the case $\textbf{(P3)}$, define a coordinate system such that $T_{x^{*}} \Sigma = \{ z = 0 \} $, the plane $T_{\theta, \lambda_{\theta}^{*}}$ coincide with $\{ x = 0 \}$, where the axis $z$ pointing to $T_{\theta, \lambda_{\theta}^{*}}$ and axis $x$ point to $\inte A$. So, use the Lemma \ref{PM2} to conclude that $w = 0$ in a neighborhood of $(0,0)$, i.e., $\widetilde{\Sigma}_{\theta, \lambda_{\theta}^{*}} = \Sigma\!\setminus\!\Sigma_{\theta, \lambda_{\theta}^{*}}$ in a neighborhood of $x^{*}$. 

\indent In $\textbf{(P4)}$, define a coordinate system such that $x^{*}$ is the origin of a coordinate system $(x, y, z)$, $T_{x^{*}} \partial\Sigma = \{ x = z = 0 \}$ and $T_{x^{*}} \Sigma = \{ z = 0 \} $, where the axis $z$ points into $\Sigma$ and axis $x$ pointing to $\widetilde{\Sigma}_{\theta, \lambda_{\theta}^{*}}$. So, use the Lemma \ref{PM3} and the capillarity of $\Sigma$ for conclude that $w = 0$ in a neighborhood of $(0,0)$, i.e., $\widetilde{\Sigma}_{\theta, \lambda_{\theta}^{*}} = \Sigma\!\setminus\!\Sigma_{\theta, \lambda_{\theta}^{*}}$ in a neighborhood of $x^{*}$.

\indent Using the unique continuation we concluded that $\widetilde{\Sigma}_{\theta, \lambda_{\theta}^{*}} = \Sigma\!\setminus\!\Sigma_{\theta, \lambda_{\theta}^{*}}$. 

\textbf{Affirmation:} $\lambda_{\theta}^{*} = 0$, $\ \forall \ \theta \in [0,\pi)$.  

\indent Indeed, suppose absurdly, that $\lambda_{\theta}^{*} < 0$. 

\indent Hereafter, as $\partial \Sigma_{\theta, \lambda_{\theta}^{*}} \subset \mathbb{S}^{2}$, then $\Sigma$ is a surface whose boundary satisfy

\begin{align}
\partial \Sigma = \partial \Sigma_{\theta, \lambda_{\theta}^{*}} \cup \partial \tilde{\Sigma}_{\theta, \lambda_{\theta}^{*}},
\end{align}

\noindent i.e., $\partial\Sigma$ does not contained in $\mathbb{S}^{2}$. Contradiction, because $\phi$ is admissible! Then, the affirmation is true.

\indent Finally, as $\lambda_{\theta}^{*} = 0$ and $\tilde{\Sigma}_{\theta, 0} = \Sigma\setminus \Sigma_{\theta, 0}$, $\forall \ \theta \in [0,\pi)$, because $\theta$ was taken arbitrarily, if $\Pi$ is a plane parallel to $\Pi_3$, the straight line $r_{\theta} := \Pi \cap T_{\theta, 0}$ intersects orthogonally $\Sigma \cap \Pi$, for all $\theta \in [0, \pi)$. Besides that, as $\lambda_{\theta}^{*} = 0$, $\forall \ \theta \in [0,\pi)$, all these straight lines intersects each other at point $p_0 \in \Pi \cap \textrm{axis}\,x_3$, $\forall \ \theta \in [0,\pi)$, i.e., $\Sigma \cap \Pi$ is a circle.

\indent Therefore, as $\theta$ was taken arbitrarily, $\Sigma$ is symmetrical rotationally.
\begin{flushright}
$\blacksquare$
\end{flushright}

\begin{figure}[ht]  
\centering
\includegraphics[scale=0.3]{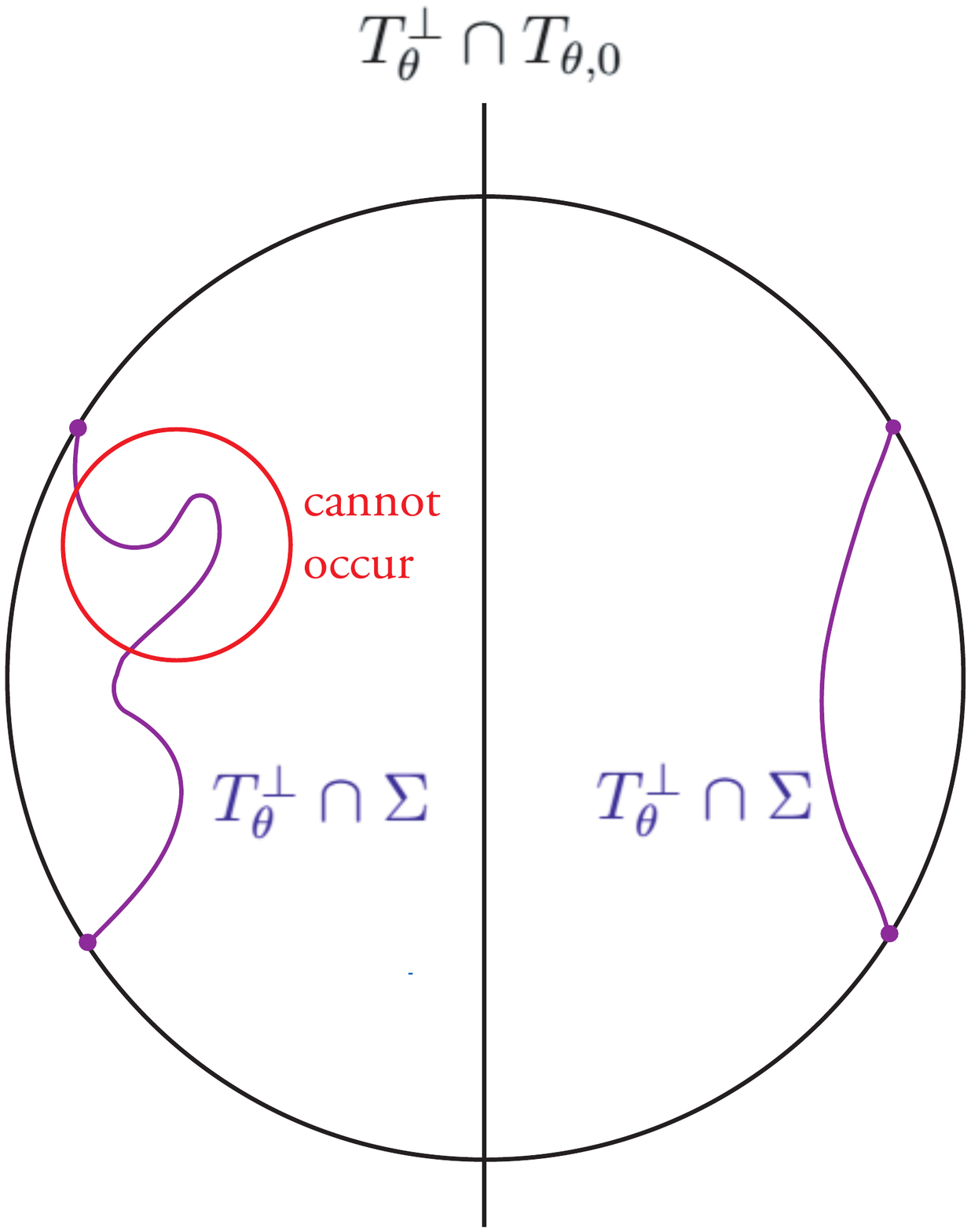}
\caption{As $\textbf{(P3)}$ does not occur for $\lambda < 0$, the circled part doesn't occur either.}
\label{secao_de_Sigma_3}
\end{figure}

\begin{observation}
Follows from above affirmation that, for example, $\textbf{(P3)}$ does not occur for $\lambda < 0$. So, the curve defined by intersection $T_{\theta}^{\perp} \cap \Sigma$, where $T_{\theta}^{\perp}$ is the plane containing the origin and orthogonal to $\Pi_{3}$ and $\mathcal{T}_{\theta}$, can be represented, globally, by the graph of a smooth function $f(z)$, where $z \in I \subset T_{\theta}^{\perp} \cap T_{\theta, 0}$, see Figure \ref{secao_de_Sigma_3}. Thus, not exist the possibility of a boundary point of $\widetilde{\Sigma}_{\lambda}$ intersects $\Sigma\!\setminus\!\Sigma_{\lambda}$, i.e., we could considered the concept \textit{touching} from the start.
\end{observation}

In Theorem \ref{Teo_McGrath}, McGrath assume that an embedded minimal free boundary annulus, $\Sigma \subset \mathbb{B}^n$, is $G$-invariant, to prove that $\Sigma$ is the critical catenoid. From Theorem \ref{Resultado Principal}, we improved the result of McGrath \cite{McGrath}, to $n=2$, because we assume only that $\partial\Sigma$ is $G$-invariant.

\begin{corollary} \label{McGrath melhorado}
Let $\Sigma^2 \subset \mathbb{B}^3$ be an embedded minimal free boundary annulus. If $\partial\Sigma$ is $G$-invariant, then $\Sigma$ is the critical catenoid.
\end{corollary}

\noindent \textbf{Proof of Corollary \ref{McGrath melhorado}:} Follows directly of proof from Theorem (\ref{Resultado Principal}) and of fact that the critical catenoid is the only minimal surface rotationally symmetric free boundary in $\mathbb{B}^3$.
\begin{flushright}
$\blacksquare$
\end{flushright}

\indent With this methodology, we also get a new demonstration for

\begin{theorem}[Pyo]  \label{Pyo_diferente}
Let $\Sigma^2$ be an embedded minimal surface in $\mathbb{R}^3$ with two boundary components and let $\Gamma$ be one component of $\partial\Sigma$. If $\Gamma$ is a circle and $\Sigma$ meets a plane along $\Gamma$ at a constant angle, then $\Sigma$ is part of the catenoid.
\end{theorem}

\noindent \textbf{Proof of Theorem \ref{Pyo_diferente}:} Let $\Pi$ the plane that contain $\Gamma$ and $\mathcal{T}_{\theta}$ be a family of parallel planes with each other and orthogonal to $\Pi$. 

\indent As $\Gamma$ is a circle, during the moving plane process, for some value of the parameter $\lambda$, $\widetilde{\Sigma}_{\theta, \lambda}$ definitely extrapolates $\Sigma\!\setminus\!\widetilde{\Sigma}_{\theta, \lambda}$, of some of the forms $(P1) \sim (P4)$. As we have no information about $\Gamma'$, another component connected of $\partial\Sigma$, we cannot say anything about the occurrence of the cases $(P2)$ and $(P4)$. If the cases $(P1)$ or $(P3)$ occur for some $\lambda_{\theta}^{*} \leq 0$, we have by ARM that
\begin{align}
\widetilde{\Sigma}_{\theta, \lambda_{\theta}^{*}} \textrm{ coincide to } \Sigma \!\setminus\! \widetilde{\Sigma}_{\theta, \lambda_{\theta}^{*}}
\end{align}

\noindent and
\begin{align}
\Gamma = \Gamma_{\theta, \lambda_{\theta}^{*}} \cup \tilde{\Gamma}_{\theta, \lambda_{\theta}^{*}},
\end{align}

\noindent But since $\Gamma$ is a circle, it follows of the circular geometry of $\Gamma$, which $\lambda_{\theta}^{*} = 0$. 

\indent Consider $r$ the orthogonal straight line to $T_{\theta, 0}$ passing through the center of $\Gamma$ and let $p$ the point given by the intersection between $\Gamma$, $r$ and $\Sigma \!\setminus\! \widetilde{\Sigma}_{\theta, 0}$. Once $\Gamma$ is a circle and $(P1)$ and $(P3)$ do not occur for $\lambda < 0$, we have $\widetilde{\Sigma}_{\theta, 0}$ stays above $\Sigma\!\setminus\!\widetilde{\Sigma}_{\theta, 0}$, relative to $\tilde{\nu}$, normal vector for $\widetilde{\Sigma}_{\theta, 0}$ at the point $p$.

\indent Consider a coordinated system such that $\{ x_3 = 0 \} = T_{x^{*}} \widetilde{\Sigma}_{\theta, 0}$. Thus, using the ARM, the unique continuation principle, the arbitrariness in choosing $\theta$, as well as in the proof of Theorem \ref{Resultado Principal}, we conclude that $\Sigma$ is the critical catenoid.
\begin{flushright}
$\blacksquare$
\end{flushright}


\begin{thebibliography}{999999}

\markboth{}{}

\bibitem{Wente} Wente, H. C. The symmetry of sessile and pendent drops, Pac. J. Math., 88 (1980),
387-397.

\bibitem{Nitsche} Nitsche, J. C. C.: Stationary partitioning of convex bodies, Arch. Rational Mech. Anal 89
(1985), 1-19.

\bibitem{Fraser_Li} Ailana Fraser and Martin Man-chun Li. Compactness of the space of embedded minimal surfaces with free boundary in three-manifolds with nonnegative Ricci curvature and convex boundary. J. Differential Geom., 96(2):183–200, 2014.

\bibitem{Brendle} Simon Brendle. Embedded minimal tori in $\mathbb{S}^3$ and the Lawson conjecture. Acta Math., 211(2):177–190, 2013.

\bibitem{Ros} Antonio Ros. A two-piece property for compact minimal surfaces in a three-sphere. Indiana Univ. Math. J.,
44(3):841–849, 1995.

\bibitem{McGrath} Peter McGrath. A characterization of the critical catenoid. arXiv preprint, arXiv:1603.04114v3.

\bibitem{Fraser_Schoen} Ailana Fraser and Richard Schoen. Sharp eigenvalue bounds and minimal surfaces in the ball. Invent. Math., 203(3):823–890, 2016.

\bibitem{Koiso} Koiso, M. Symmetry of Hypersurfaces of Constant Mean Curvature with Symmetric Boundary. Math. Z. 191,567-574 (1986).

\bibitem{Pyo} Juncheol Pyo. Minimal annuli with constant contact angle along the planar boundaries. Geom Dedicata (2010) 146: 159.

\bibitem{Serrin} James Serrin. A symmetry problem in potential theory. Arch. Rat. Mech. and Anal., 43
(1971), 304-318.

%\bibitem{Cheng_1} Cheng, S. Y.: Eigenfunctions and eigenvalues of Laplacian. pages 185–193, 1975.

%\bibitem{Cheng_2} Cheng, S. Y.: Eigenfunctions and nodal sets. Comment. Math. Helv. 51 (1976),43-55.



\end{thebibliography}
\end{document}